\UseRawInputEncoding
\documentclass{amsart}
\usepackage{amsmath,amssymb}

\newtheorem{example}{Example}

\newtheorem{conjecture}{Conjecture}

\usepackage{epsfig}
\usepackage{url}
\usepackage{amssymb}
\usepackage{amsmath}
\usepackage{amsfonts}

\def\Dbar{\leavevmode\lower.6ex\hbox to 0pt{\hskip-.23ex \accent"16\hss}D}

\def\bZ{{\mbox{\bf Z}}}

\begin{document}

\title{Some new symmetric Hadamard matrices}
\author {Dragomir {\v{Z}}. {\Dbar}okovi{\'c}}
\address{University of Waterloo, 
Department of Pure Mathematics and Institute for Quantum Computing,
Waterloo, Ontario, N2L 3G1, Canada}
\email{djokovic@uwaterloo.ca}
\date{}

\begin{abstract}
The first examples of symmetric Hadamard matrices of orders 
$4\cdot 127$ and $4\cdot 191$ are presented.
The systematic computer search for symmetric Hadamard matrices 
based on the so called propus array 
was carried out recently for all orders $4v$ with odd $v\le53$.
Hypothetically, such matrices exist for all odd $v$ and all 
propus parameter sets $(v;k_1,k_2,k_3,k_4;\lambda)$, $k_2=k_3$,
$\sum k_i=\lambda+v$, apart from the exceptional cases when all the $k_i$ are equal. In this note the search has been extended 
further to cover the cases $v=55,57,59,61,63$.
\end{abstract}

\maketitle

\section{Introduction}

We fix some notation which will be used throughout this note.
Let $X_i$, $i=1,2,3,4$, be a difference family (DF) in a finite abelian group $G$ (written additively) and let 
$$(v;k_1,k_2,k_3,k_4;\lambda)$$
be its parameter set (PS). Thus $v=|G|$, $|X_i|=k_i$ and 
$\sum k_i(k_i-1)=\lambda(v-1)$, where $|X|$ denotes the cardinality of a finite set $X$. If $\sum k_i=\lambda+v$ we say that this PS is a {\em Goethals-Seidel parameter set} (GSPS) and  that this DF is a {\em Goethals-Seidel difference family} (GSDF). If the $X_i$ form a GSDF and we replace one of the 
blocks by its set-theoretic complement in $G$, we again 
obtain a GSDF although the parameter set may change. For that 
reason we shall always assume that all the $k_i$ are 
$\le v/2$.

Each GSDF in $G$ gives a Hadamard matrix $H$ of order $4v$. For more details about this construction see e.g. 
\cite{DK:CompMethods:2019,SY:1992}. Briefly, each $X_i$ provides 
a $G$-invariant matrix $A_i$ of order $v$, and $H$ is obtained 
by plugging the $A_i$ into the well known 
{\em Goethals-Seidel array} 
\begin{equation*} 
GSA=\left[ \begin{array}{cccc}
A_1    & A_2 R & A_3 R & A_4 R \\
-A_2 R & A_1   & -RA_4 & RA_3  \\
-A_3 R & RA_4  & A_1   & -RA_2 \\
-A_4 R & -RA_3 & RA_2  & A_1  \\
\end{array} \right].
\end{equation*} 
We recall that a matrix $A=(a_{x,y})$ with indices $x,y\in G$ is {\em $G$-invariant} if $a_{x+z,y+z}=a_{x,y}$ for all 
$x,y,z\in G$. The matrix $R=(r_{x,y})$ may be defined by the 
formula $r_{x,y}=\delta_{x+y,0}$, $x,y\in G$, where $\delta$ is the Kronecker symbol.

For a subset $X$ of $G$, we say that it is {\em symmetric} if 
$-X=X$, and we say that it is {\em skew} if $G$ is a disjoint union of $X$, $-X$ and $\{0\}$. If at least one of the blocks 
$X_i$ of a GSDF is skew then, after rearranging the $X_i$ so 
to have $X_1$ skew, the Hadamard matrix $H$ will be skew Hadamard, i.e. such that $H+H^T=2I_{4v}$. (T denotes the matrix transposition, and $I_k$ is the identity matrix of order $k$.)

In order to obtain a symmetric Hadamard matrix $H$ we require that two of the blocks $X_i$ are the same and that one of the other two blocks is symmetric. 
A {\em propus parameter set} (PPS) is a GSPS having $k_i=k_j$
for some $i\ne j$. By permuting the $k_i$'s we may assume that 
$k_2=k_3$ and $k_1\ge k_4$. In that case we say that this PPS 
is {\em normalized}. Note that these conditions in general do not specify the $k_i$'s uniquely. For instance the PPSs $(5;1,2,2,1;1)$ and $(5;2,1,1,2;1)$ are both normalized but they become the same if we ignore the ordering of the $k_i$'s. 

We say that a GSDF is a {\em propus difference family} (PDF) if 
$X_i=X_j$ for some $i \ne j$ and one of the other two blocks is symmetric. 
If the $X_i$'s form a PDF then, after rearranging the blocks we may assume that $X_2=X_3$ and that $X_1$ is symmetric. Then we 
plugg the corresponding matrix blocks $A_i$ into the so called 
propus array (PA) to construct a symmetric Hadamard matrix of order $4v$. This construction is known as the {\em propus construction}. It has been first introduced in \cite{SB:2017}.
For the reader's convenience we display the propus array
\begin{equation*}
PA=\left[ \begin{array}{cccc}
-A_1  & A_2 R & A_3 R & A_4 R \\
A_3 R & RA_4  & A_1   & -RA_2 \\
A_2 R & A_1   & -RA_4 & RA_3  \\
A_4 R & -RA_3 & RA_2  & A_1  \\
\end{array} \right].
\end{equation*} 
Note that PA is obtained from GSA by multiplying the first column by $-1$ and interchanging the second and the third rows. 

From now on we assume that $G=\bZ_v$, a cyclic group of order 
$v$, and that $v$ is odd. Under this assumption, the matrix blocks $A_i$ will be circulants. All PPSs for $v\le 41$ are listed in \cite{BBDKS:2017} together with the corresponding PDF's. There was only one case of a PPS having no PDF, namely $(25;10,10,10,10;15)$. Similarly, the cases $41<v\le 51$ were handled in \cite{BDK:2018}, and the case $v=53$ in \cite{ABDK:2019}. Again there was one exceptional case, 
$(49;21,21,21,21;35)$. In the present note, for each PPS 
with $51<v\le 63$ we exhibit at least one PDF. For more information on the exceptional cases see \cite{BDK:2018}.

\section{Symmetric Hadamard matrices of orders 508 and 764}
\label{glavni}

The {\em symmetry symbol} (abc) written immediately after a PPS shows the symmetry types of the three blocks $X_1$, $X_2$ and $X_4$. More precisely, the letter $s$ means that we require that the corresponding block be symmetric, the letter $k$ is used if we require that block to be skew, and the symbol ${*}$ is used otherwise. In particular $a=s$ means that we require $X_1$ to be  symmetric, $a=k$ means that we require $X_1$ to be skew, and $a={*}$ means that no symmetry condition is imposed on $X_1$. 

The group of units $\bZ^*_v$ acts on $\bZ_v$ by multiplication.
It may happen that there is a nontrivial subgroup $H$ of $\bZ^*$
such that some block $X_i$ of a PDF is a union of orbits of $H$.
In such case we may specify $X_i$ by writing it as $HY_i$, 
where $Y_i$ is a set of representatives of the $H$-orbits contained in $X_i$.

For $v=127$ we give five nonequivalent PDFs and for $v=191$ only one.

\begin{eqnarray*}
&& \qquad \qquad v=127,~ H=\{1,19,107\} \\
&& (127;57,61,61,55;107) ~ ({**}s) \\
X_1 &=& H\{4,5,6,9,12,15,23,24,30,33,36,39,45,52,58,59,60,64,
           66\} \\
X_2 &=& H\{0,4,5,6,13,15,17,26,30,32,40,46,51,53,58,59,60,64,
           65,66,72\} \\
X_4 &=& H\{0,2,4,8,9,12,15,23,24,26,30,33,40,46,51,52,53,65,
           71\} \\
\\
&& (127;60,60,60,54;107) ~ (s{**}) \\
X_1 &=& H\{1,5,6,11,13,15,16,17,20,23,24,29,32,45,46,52,58,66,
           71,72\} \\
X_2 &=& H\{2,4,5,11,12,15,16,18,22,23,29,33,36,39,46,51,52,53,
           60,71\} \\
X_4 &=& H\{6,8,17,20,22,23,30,33,36,39,45,51,58,59,60,64,
           66,71\} \\
\\
&& (127;60,60,60,54;107) ~ ({**}s) \\
X_1 &=& H\{1,2,3,4,5,6,9,11,12,13,15,17,23,24,32,33,39,46,
           64,65\} \\
X_2 &=& H\{2,5,9,10,12,13,15,16,17,29,33,36,39,40,45,51,53,
           58,60,66\} \\
X_4 &=& H\{1,5,6,10,11,13,16,17,20,23,30,32,45,58,64,65,66,
           71\} \\
\\
&& (127;58,60,60,55;106) ~ ({**}s) \\
X_1 &=& H\{0,2,3,4,5,12,13,16,17,18,20,22,29,30,46,51,53,58,
           59,71\} \\
X_2 &=& H\{8,9,10,16,20,22,23,24,26,29,32,36,45,46,51,52,59,60,
           65,78\} \\
X_4 &=& H\{0,3,5,10,11,17,18,22,24,29,32,39,45,52,58,59,60,
           64,72\} \\
\\
&& (127;60,57,57,58;105) ~ (s{**}) \\
X_1 &=& H\{2,4,6,10,12,13,15,17,23,24,26,36,40,46,51,52,58,64,
           71,78\} \\
X_2 &=& H\{1,2,3,4,10,16,17,18,20,23,29,30,45,51,52,58,64,
           66,72\} \\
X_4 &=& H\{0,2,5,8,9,10,11,13,15,17,18,30,39,40,46,53,58,60,
           66,78\} \\
\end{eqnarray*}

\begin{eqnarray*}
&& \qquad \qquad v=191,~ H=\{1,39,49,109,184\} \\
&& (191;91,90,90,85;165) ~ (s{**}) \\
X_1 &=& H\{0,1,3,4,7,9,16,17,18,21,22,28,31,36,57,61,62,
           68,112\} \\
X_2 &=& H\{1,4,14,16,18,19,22,23,28,29,31,32,34,36,38,61,
           62,68\} \\
X_4 &=& H\{1,2,9,11,12,17,18,22,28,29,31,32,38,41,56,
           61,66\} \\
\end{eqnarray*}

\section{Small orders of symmetric Hadamard matrices}

There are several known infinite series of PDFs 
\cite{DDK:SpecMatC:2015,SB:2017}. We shall use only two of them.
The first one is essentially the Turyn series \cite{Turyn:1972} 
with $v=(q+1)/2$, $q$ a prime power $\equiv 1 \pmod{4}$, and all
four blocks $X_i$ symmetric. The second one is essentially 
the series constructed in \cite{XXSW:2005} (see also \cite{DDK:SpecMatC:2015}) to which we refer as the XXSW-series. In this case $v=(q+1)/4$, $q$ a prime power $\equiv 3 \pmod{4}$, and we may arrange the blocks so that $X_1$ is skew, 
$X_2=X_3$ and $X_4$ is symmetric.

In the handbook \cite{CK-Had:2007} published in 2007 it is indicated (see Table 1.52, p.277) that, for odd $v<200$, no symmetric Hadamard matrices of order $4v$ are known for 
\begin{eqnarray*} 
v&=& 23,29,39,43,47,59,65,67,73,81,89,93,101,103,107,109,113, \\
  && 119,127,133,149,151,153,163,167,179,183,189,191,193.
\end{eqnarray*}

The cases $v=23$ and $v=81$ should not have been included. For the case $v=23$ see \cite{DDK:SpecMatC:2015}. For $v=81$ note that symmetric Hadamard matrices of orders $4\cdot 9^k$, 
$k\ge 1$ integer, were constructed by Turyn \cite{Turyn:JCT:1984} back in 1984. Moreover, the Bush-type Hadamard matrix of order $4\cdot 81=324$ constructed in 2001 \cite{JKT:2001} is also symmetric.

The propus construction has been used in several recent papers \cite{ABDK:2019,BBDKS:2017,BD:2018,BDK:2018,DDK:SpecMatC:2015,SB:2017} to construct symmetric Hadamard matrices of new orders. By taking into account these results and those from Sect. \ref{glavni}, the above list of undecided cases reduces to
\begin{eqnarray*} 
v&=& 65,89,93,101,107,119,133,149,153,163,167,179,183,189,193.
\end{eqnarray*}

\section{List of PPSs and PDFs for odd $v$, $53 < v \le 63$}

The following conventions and notation will be used in the listings below. We have $\bZ_v=\{0,1,2,...,v-1\}$ and recall that $v$ is odd. Let $X \subseteq \bZ_v$ and $k=|X|$. Define 
$X' = X \cap \{1,2,...,(v-1)/2\}$. In particular 
$\bZ'_v = \{1,2,...,(v-1)/2\}$.

If $X$ is skew then $k=(v-1)/2$ and 
$$X = X' \cup (-(\bZ'_v \setminus X')).$$

If $X$ is symmetric then  
$$ 
X = 
\begin{cases} X' \cup (-X'), & {\rm for} ~ k~ {\rm even}; \cr 
     \{0\} \cup X' \cup (-X'), & {\rm for} ~ k ~ {\rm odd}. \cr \end{cases}
$$

Hence, a skew $X$ can be recovered uniquely from $X'$. This is also true for symmetric $X$ provided we know the parity of $k$.

For a PDF $X_i$, $i=1,2,3,4$, with normalized PPS 
$(v;k_1,k_2,k_3,k_4;\lambda)$ we always assume that $X_2=X_3$. Thus it suffices to specify only the blocks $X_1$, $X_2$ and 
$X_4$. We say that a PPS is {\em exceptional} if all the $k_i$ are equal. The following conjecture is implicit in 
\cite{ABDK:2019,BBDKS:2017,BDK:2018}. It has been verified there for odd $v \le 53$.

\begin{conjecture}
For each normalized and non-exceptional PPS 
$(v;k_1,k_2,k_3,k_4;\lambda)$ there exist PDFs with symmetry 
symbols $(s{**})$ and $({**}s)$.
\end{conjecture}

The list below shows that the conjecture is true also for 
$v=55,57,59,61,63$. 

If a block $X_i$ is symmetric or skew, in order to save space we record only $X'_i$. As the $k_i$ are specified by the PPS, $X_i$ can be recovered uniquely from $X'_i$. 

\begin{example}
For the first PDF below, the symmetry symbol $(s{**})$ shows that $X_1$ must be symmetric. As 
$X'_1=\{5,6,7,9,10,13,15,16,19,21,23,25,27\}$ 
we have $-X'_1=\{28,30,32,34,36,39,40,42,45,46,48,49,50\}$. 
As $k_1=27$ is odd we have $X_1=\{0\} \cup X'_1 \cup (-X'_1)$.
\end{example}

\newpage
\begin{eqnarray*}
&& \qquad \qquad v=55 \\
&& (55;27,25,25,21;43) ~ (s{**}) \\
X'_1 &=& \{5,6,7,9,10,13,15,16,19,21,23,25,27\} \\
X_2 &=& \{0,1,2,3,4,6,9,10,14,17,19,24,26,29,30,34,37,38,39,40,
41,47,48, \\
&& ~ 52,53\} \\
X_4 &=& \{0,3,4,5,10,11,12,14,16,17,18,19,21,22,24,30,36,43,44,
46,47\} \\
\\
&& (55;27,25,25,21;43) ~ ({**}s) \\
X_1 &=& \{0,4,5,6,7,9,10,12,15,20,21,24,25,26,28,32,33,34,38,39,
41,44,45, \\
&& ~ 51,52,53,54\} \\
X_2 &=& \{0,3,5,7,8,9,10,19,23,25,28,31,32,33,34,36,37,40,41,43,
44,45,47, \\
&& ~ 48,53\} \\
X'_4 &=& \{2,4,6,9,12,13,18,19,20,23\} \\
\\
&& (55;27,24,24,22;42) ~ (s{**}) \\
X'_1 &=& \{1,4,6,7,9,11,14,15,18,20,23,24,27\} \\
X_2 &=& \{0,4,6,7,12,14,15,16,18,23,25,26,30,31,33,34,35,36,37,
40,41,46, \\
&& ~ 47,53\} \\
X_4 &=& \{0,1,4,6,16,17,18,19,20,21,22,23,24,27,29,31,33,36,42,
43,44,46\} \\
\\
&& (55;27,24,24,22;42) ~ ({**}s) \\
X_1 &=& \{0,1,5,8,10,12,15,16,17,24,25,26,29,30,34,37,39,40,41,
42,44,47, \\
&& 48,50,52,53,54\} \\
X_2 &=& \{0,10,13,14,15,16,17,20,21,22,24,26,30,31,33,34,35,41,
43,44,47, \\
&& ~ 49,50,53\} \\
X'_4 &=& \{1,6,11,12,13,16,19,20,22,24,27\} \\
\\
&& (55;26,23,23,24;41) ~ (s{**}) \\
X'_1 &=& \{1,2,5,8,10,11,13,14,15,19,21,23,27\} \\
X_2 &=& \{0,2,3,4,5,6,7,10,11,14,18,24,25,30,35,37,39,40,41,42,
50,51,52\} \\
X_4 &=& \{0,2,5,6,8,14,16,19,20,21,22,24,25,28,31,32,33,37,38,40,
45,47, \\
&& ~ 49,52\} \\
\\
&& (55;26,23,23,24;41) ~ ({**}s) \\
X_1 &=& \{0,2,7,11,14,15,17,18,23,24,25,28,29,30,31,36,37,38,39,
42,44,45,\\
&& ~ 46,48,50,53\} \\
X_2 &=& \{0,3,8,12,13,14,18,19,23,26,33,34,36,38,42,45,46,47,48,
49,50, \\
&& ~ 51,52) \\
X'_4 &=& \{1,2,7,10,11,13,17,19,21,24,26,27\} \\
\end{eqnarray*}

\begin{eqnarray*}
&& (55;24,27,27,21;44) ~ (s{**}) \\
X'_1 &=& \{6,8,10,13,15,16,18,19,20,22,25,26\} \\
X_2 &=& \{0,2,4,5,8,13,14,16,17,19,25,26,27,32,33,34,37,38,39, 40,41, \\
&& ~ 42,44,49,50,53,54\} \\
X_4 &=& \{0,4,5,8,11,12,13,16,18,20,22,24,31,33,34,36,37,41,42, 44,51\} \\
\\
&& (55;24,27,27,21;44) ~ ({**}s) \\
X_1 &=& \{0,3,11,12,13,14,16,17,23,24,26,27,29,30,35,38,39,40,
41,42,47,48, \\
&& 49,50\} \\
X_2 &=& \{0,1,3,4,5,6,8,9,14,16,18,21,22,27,28,32,35,36,37,39,43,
47,49, \\
&& ~ 51,52,53,54\} \\
X'_4 &=& \{3,5,9,10,14,16,17,20,23,25\} \\
\\
&& (55;24,25,25,22;41) ~ (s{**}) \\
X'_1 &=& \{1,3,5,6,7,11,14,15,16,18,21,25\} \\
X_2 &=& \{0,9,11,12,14,17,20,22,23,24,25,26,27,30,31,33,37,38,42,46,
47,48,\\
&& ~ 49,52,54\} \\
X_4 &=& \{0,4,5,7,8,11,12,16,18,19,21,24,25,26,28,34,35,39,41,
43,53,54\} \\
\\
&& (55;24,25,25,22;41) ~ ({**}s) \\
X_1 &=& \{0,2,3,5,13,14,16,17,21,22,23,26,29,32,37,38,42,43,44,
45,46,48, \\
&& ~ 49,51\} \\
X_2 &=& \{0,1,2,3,4,10,11,15,17,18,19,21,22,24,25,27,28,30,32,
34,35,39, \\
&& ~ 40,44,46\} \\
X'_4 &=& \{2,3,5,7,10,11,15,21,23,25,26\} \\
\\
&& (55;23,26,26,22;42) ~ (sss), ~  {\rm Turyn ~ series} \\
X'_1 &=& \{6,7,10,11,15,17,18,19,21,24,26\} \\
X'_2 &=& \{1,2,4,8,14,16,17,18,19,23,24,25,27\} \\
X'_4 &=& \{6,7,10,11,15,17,18,19,21,24,26\} \\
\\
&& \qquad \qquad v=57, ~ H=\{1,7,49\} \\
&& (57;28,28,28,21;48) ~ (s{**}) \\
X'_1 &=& \{1,2,6,8,10,11,14,16,17,18,21,22,25,27\} \\
X_2 &=& \{0,2,3,5,6,7,11,12,13,15,18,23,24,25,26,27,28,29,30,32,
34,39,40, \\
&& 42,46,49,50,55\} \\
X_4 &=& \{0,1,2,3,5,6,8,9,12,13,14,15,19,23,30,31,32,39,45,48,51\} \\
\end{eqnarray*}

\begin{eqnarray*}
&& (57;28,28,28,21;48) ~ (k{*}s), ~ {\rm XXSW - series} \\
X'_1 &=& \{2,4,12,13,15,21,23,24,25,27,28\} \\
X_2 &=& \{1,3,5,8,9,12,15,20,23,24,26,27,29,31,32,33,35,36,37,
41,42,45,49, \\
&& ~ 50,51,52,53,55\} \\
X'_4 &=& \{1,4,6,13,14,15,19,20,21,26\} \\
\\
&&(57;27,26,26,22;44) ~ (s{**}), ~ 
{\rm all} ~ X_i ~ {\rm are} ~ H{\rm {}-invariant} \\
X'_1 &=& \{3,6,9,10,11,13,15,19,20,21,23,24,26\} \\
X_2 &=& \{4,6,8,9,15,16,19,23,25,28,30,31,37,38,39,42,43,44,45,
46,47,48, \\
&& ~ 50,51,55,56\} \\
X_4 &=& \{2,4,6,9,10,11,13,14,20,22,25,26,28,30,34,38,39,40,41,
42,45,52\} \\
\\
&& (57;27,26,26,22;44) ~ ({**}s) \\
X_1 &=& \{0,8,9,10,12,13,15,18,19,20,25,30,33,34,37,40,41,43,47,
48,49,51, \\
&& ~ 52,53,54,55,56\} \\
X_2 &=& \{0,4,8,9,11,16,17,18,22,24,25,26,27,28,30,31,33,34,37,
38,42,43, \\
&& ~ 51,53,54,55\} \\
X'_4 &=& \{2,4,7,9,11,12,14,17,21,23,28\} \\
\\
&& (57;27,25,25,23;43) ~ (s{**}), ~ 
{\rm all} ~ X_i ~ {\rm are} ~ H{\rm {}-invariant} \\
X'_1 &=& \{3,6,9,10,11,13,15,19,20,21,23,24,26\} \\
X_2 &=& \{2,3,4,14,16,21,22,24,25,28,29,30,32,33,36,38,39,40,41,
43,45,52, \\
&& 53,54,55\} \\
X_4 &=& \{1,2,5,7,10,13,14,17,19,24,29,30,32,34,35,36,38,39,41,
45,49, \\
&& ~ 53,54\} \\
\\
&& (57;27,25,25,23;43) ~ ({**}s) \\
X_1 &=& \{0,1,7,11,15,16,17,19,20,21,24,25,26,28,30,35,36,37,38,
40,44,46, \\
&& ~ 47,48,49,51,54\} \\
X_2 &=& \{0,1,2,3,7,10,14,17,19,20,23,26,31,32,34,36,37,38,41,
42,43,44, \\
&& ~ 45,46,49\} \\
X'_4 &=& \{2,4,6,9,10,11,15,16,18,23,26\} \\
\\
&& (57;25,25,25,24;42) ~ (sss), ~ {\rm Turyn ~ series} \\
X'_1 &=& \{2,3,8,9,10,18,20,22,23,24,26,27\} \\
X'_2 &=& \{6,7,9,10,14,16,19,21,24,25,27,28\} \\
X'_4 &=& \{2,3,8,9,10,18,20,22,23,24,26,27\} \\
\end{eqnarray*}

The third and the sixth PDF below are taken from \cite{ABDK:2019}.
\begin{eqnarray*}             
&& \qquad \qquad v=59 \\
&& (59;28,29,29,22;49) ~ (s{**}) \\
X'_1 &=& \{1,3,5,8,10,11,13,15,16,20,21,22,26,29\} \\
X_2 &=& \{0,1,3,7,8,9,10,12,13,15,16,19,21,22,25,27,29,34,35,36,
37,38, \\
&& ~ 39,40,44,51,54,55,58\} \\
X_4 &=& \{0,3,4,5,7,8,14,15,16,17,19,22,24,27,28,33,39,49,53,54,
55,56\} \\
\\
&& (59;28,29,29,22;49) ~ ({**}s) \\
X_1 &=& \{0,5,6,7,8,9,10,12,17,18,20,25,26,27,34,35,39,42,44,45,
47,48,49, \\
&& ~ 50,51,54,55,58\} \\
X_2 &=& \{0,1,2,3,4,5,6,9,11,13,14,18,20,24,25,26,29,31,32,35,
37,44,45,47, \\
&& 48,51,55,56,58\} \\
X'_4 &=& \{1,5,7,11,15,16,18,20,21,23,29\} \\
\\
&& (59;27,25,25,26;44) ~ (s{**}) \\
X'_1 &=& \{2,4,7,8,12,13,15,16,17,18,20,23,29\} \\
X_2 &=& \{1,2,4,5,12,13,17,19,20,21,22,23,26,27,31,35,37,38,40,
44,47,49, \\
&& ~ 50,55,57\} \\
X_4 &=& \{3,7,12,13,14,16,18,19,20,22,23,24,25,26,31,32,33,34,36,
38,43,45, \\
&& ~ 46,50,51,53\} \\
\\
&& (59;27,25,25,26;44) ~ ({**}s) \\
X_1 &=& \{0,1,3,5,6,7,8,9,12,15,18,20,28,29,31,33,34,35,38,42,
44,47,48,49, \\
&& ~ 55,56,58\} \\
X_2 &=& \{0,3,4,5,7,10,16,21,22,24,25,26,28,29,32,33,34,38,39,
40,41,43,48, \\
&& ~ 49,52\} \\
X'_4 &=& \{1,3,5,7,10,12,13,15,18,19,20,27,28\} \\
\\
&& (59;26,28,28,23;46) ~ (s{**}) \\
X'_1 &=& \{4,6,10,12,13,15,17,21,22,24,25,27,29\} \\
X_2 &=& \{0,1,4,5,6,7,10,11,12,13,14,17,21,22,23,27,33,34,36,37,
39,41,42, \\
&& ~ 45,52,55,56,57\} \\
X_4 &=& \{0,7,9,12,13,14,22,25,26,27,28,31,33,35,39,40,46,47,49,
51,55, \\
&& ~ 57,58\} \\
\end{eqnarray*}

\begin{eqnarray*}               
&& (59;26,28,28,23;46) ~ ({**}s) \\
X_1 &=& \{2,3,10,12,13,14,16,18,19,26,28,29,36,38,39,40,42,44,
46,47,49, \\
&& ~ 50,53,54,55,57\} \\
X_2 &=& \{4,5,7,11,12,16,17,24,25,26,27,28,29,33,34,37,39,40,42,
43,44,45, \\
&& ~ 47,49,51,53,56,58\} \\
X'_4 &=& \{1,4,5,7,8,11,14,20,25,28,29\} \\
\\
&& \qquad \qquad v=61,~H_1=\{1,13,47\},~H_2=\{1,9,20,34,58\} \\ && (61;30,29,29,23;50) ~ (s{**}) \\
X'_1 &=& \{2,3,8,10,11,13,14,16,18,20,22,23,24,28,30\} \\
X_2 &=& \{0,1,2,4,9,10,11,13,18,19,22,24,26,27,37,38,39,40,42,
43,44,48,49, \\
&& ~ 51,52,55,56,58,59\} \\
X_4 &=& \{0,6,7,8,11,12,14,16,17,23,26,30,31,33,34,35,37,38,40,
43,45,49,50\} \\
\\
&& (61;30,29,29,23;50) ~ ({**}s) \\
X_1 &=& \{1,2,5,6,7,8,11,12,14,15,19,20,21,22,24,27,28,30,31,32,
38,39,40, \\
&& ~ 41,43,45,46,55,58,60\} \\
X_2 &=& \{0,2,4,8,9,10,13,14,17,18,20,23,24,25,26,27,32,37,38,
44,45,47, \\
&& ~ 49,53,55,56,58,59,60\} \\
X'_4 &=& \{2,10,13,15,17,18,19,21,22,26,29\} \\
\\
&& (61;30,26,26,26;47) ~ (s{**}), ~ 
{\rm all} ~ X_i ~ {\rm are} ~ H_2{\rm {}-invariant} \\
X'_1 &=& \{1,3,4,5,9,12,13,14,15,16,19,20,22,25,27\} \\
X_2 &=& \{0,1,6,8,9,11,12,20,21,25,26,28,30,32,34,37,38,42,43,
44,47,51, \\
&& ~ 54,57,58,59\} \\
X_4 &=& \{0,3,5,6,21,23,24,26,27,30,32,33,39,41,43,44,45,46,48,
50,51,52, \\
&& ~ 53,54,59,60\} \\
\\
&& (61;30,26,26,26;47) ~ ({**}s) \\
X_1 &=& \{0,2,4,5,6,8,9,10,12,18,19,21,22,23,25,26,27,28,30,32,
34,36,37, \\
&& ~ 42,49,55,56,57,58,59\} \\
X_2 &=& \{0,2,5,8,9,10,15,18,24,27,28,31,33,35,38,39,44,45,46,
47,49,50,51, \\
&& ~ 52,59,60\} \\
X'_4 &=& \{2,3,6,9,11,14,18,19,21,23,24,25,29\} \\
\\
&& (61;30,25,25,30;49) ~ (sss), ~ {\rm Turyn ~ series} \\
X'_1 &=& \{1,2,6,8,9,12,13,14,15,16,17,19,24,25,28\} \\
X'_2 &=& \{1,5,6,8,10,11,12,14,20,24,27,29\} \\
X'_4 &=& \{3,4,5,7,10,11,18,20,21,22,23,26,27,29,30\} \\
\end{eqnarray*}

\begin{eqnarray*}               
&& (61;30,25,25,30;49) ~ (k{*}s), ~ {\rm XXSW ~ series}, \\
X'_1 &=& \{3,4,6,13,14,15,16,19,21,22,23,24,26,27,30\} \\
X_2 &=& \{5,6,8,12,14,15,17,20,31,32,33,36,40,44,45,46,48,49,51,
53,54,55, \\
&& ~ 56,59,60\} \\
X'_4 &=& \{1,3,5,9,10,13,15,16,17,20,22,26,27,29,30\} \\
\\
&& (61;28,28,28,24;47) ~ (s{**}) \\
X'_1 &=& \{1,4,6,7,8,10,11,14,19,20,22,26,28,30\} \\
X_2 &=& \{0,1,6,10,12,13,17,21,22,23,26,27,29,32,34,35,36,39,40,
42,43,50, \\
&& ~ 52,54,57,58,59,60\} \\
X_4 &=& \{0,5,6,8,18,23,24,28,32,34,37,42,43,44,48,49,50,51,52,
54,55,56, \\
&& ~ 57,58\} \\
\\
&& (61;28,28,28,24;47) ~ ({**}s), ~ 
{\rm all} ~ X_i ~ {\rm are} ~ H_1{\rm {}-invariant} \\
X_1 &=& \{0,1,3,4,5,9,12,13,14,15,18,19,27,32,34,39,40,46,47,48,
49,50,51, \\
&& ~ 52,53,56,57,60\} \\
X_2 &=& \{0,1,2,4,5,7,11,13,14,21,22,24,26,29,30,31,33,36,37,41,
42,45,47, \\
&& ~ 48,52,54,58,60\} \\
X'_4 &=& \{1,3,11,13,14,16,19,20,21,22,25,29\} \\
\\
&& (61;28,27,27,25;46) ~ (s{**}) \\
X'_1 &=& \{1,2,3,4,5,7,9,11,17,18,20,24,27,29\} \\
X_2 &=& \{0,4,6,7,11,12,15,18,19,21,26,30,31,32,33,35,36,41,43,
44,45,48, \\
&& ~ 49,51,52,53,54\} \\
X_4 &=& \{0,4,10,11,13,16,17,21,22,27,30,32,37,38,40,42,43,44,46,47,
51, \\
&& ~ 53,54,55,56\} \\
\\
&& (61;28,27,27,25;46) ~ ({**}s), ~ 
{\rm all} ~ X_i ~ {\rm are} ~ H_1{\rm {}-invariant} \\
X_1 &=& \{0,7,11,18,21,22,23,24,28,29,30,31,32,35,36,37,40,41,
42,44,45,50, \\
&& ~ 51,53,54,55,58,59\} \\
X_2 &=& \{1,3,7,8,9,10,13,19,23,24,27,28,30,31,35,37,39,43,44,
46,47,49,54, \\
&& ~ 55,56,57,59\} \\
X'_4 &=& \{1,3,8,10,13,14,16,18,19,20,22,25\} \\
\end{eqnarray*}

\begin{eqnarray*}               
&& (61;25,30,30,25;49) ~ (s{**}), ~ 
{\rm all} ~ X_i ~ {\rm are} ~ H_1{\rm {}-invariant} \\
X'_1 &=& \{1,3,8,10,13,14,16,18,19,20,22,25\} \\
X_2 &=& \{2,6,7,14,17,18,22,23,24,26,27,28,30,33,35,36,38,41,
42,44,45,46, \\
&& ~ 48,49,51,53,55,58,59,60\} \\
X_4 &=& \{0,1,2,3,6,13,14,17,19,26,27,31,32,33,37,38,39,40,46,
47,48,49,50, \\
&& ~ 54,60\} \\
\\
&& \qquad \qquad v=63, ~H_1=\{1,4,16\},~H_2=\{1,25,58\} \\
&& (63;31,26,26,30;50) ~ (sss), ~ {\rm Turyn ~ series},
{\rm all} ~ X_i ~ {\rm are} ~ H_2{\rm {}-invariant} \\
X'_1 &=& \{1,3,4,5,7,12,14,15,19,20,25,26,28,29,31\} \\
X'_2 &=& \{7,8,9,11,14,18,19,21,23,27,28,29,31\} \\
X'_4 &=& \{1,3,4,5,7,12,14,15,19,20,25,26,28,29,31\} \\
\\
&& (63;30,30,30,24;51) ~ (s{**}), ~ 
{\rm all} ~ X_i ~ {\rm are} ~ H_1{\rm {}-invariant} \\
X'_1 &=& \{1,2,4,8,9,11,13,16,18,19,22,25,26,27,31\} \\
X_2 &=& \{7,9,11,13,14,15,18,19,22,25,26,27,28,30,35,36,37,38,
39,41,44,45, \\
&& 49,50,51,52,54,56,57,60\} \\
X_4 &=& \{5,9,13,15,17,18,19,20,22,23,25,26,29,31,36,37,38,41,
51,52,53,55, \\
&& ~ 60,61\} \\
\\
&& (63;30,30,30,24;51) ~ ({**}s), ~ 
{\rm all} ~ X_i ~ {\rm are} ~ H_1{\rm {}-invariant} \\
X_1 &=& \{3,5,6,9,10,12,13,14,17,18,19,20,23,24,26,29,30,33,34,
35,36,38, \\
&& ~ 39,40,41,48,52,53,56,57\} \\
X_2 &=& \{3,5,7,9,10,12,13,15,17,18,19,20,23,26,27,28,29,34,36,
38,40,41, \\
&& ~ 45,48,49,51,52,53,54,60\} \\
X'_4 &=& \{5,7,9,10,14,17,18,20,23,27,28,29\} \\
\\
&& (63;30,27,27,27;48) ~ (s{**}), ~ 
{\rm all} ~ X_i ~ {\rm are} ~ H_2{\rm {}-invariant} \\
X'_1 &=& \{1,5,8,9,11,16,17,18,19,22,23,25,27,29,31\} \\
X_2 &=& \{3,4,7,12,15,16,17,20,22,26,27,28,29,32,37,41,43,44,
45,46,47,48, \\
&& ~ 49,51,54,59,60\} \\
X_4 &=& \{4,6,9,10,13,17,18,19,24,27,29,31,32,33,34,36,37,40,
41,43,44,45, \\
&& ~ 47,52,54,55,61\} \\
\end{eqnarray*}

\begin{eqnarray*}
&& (63;30,27,27,27;48) ~ ({**}s), ~  
{\rm all} ~ X_i ~ {\rm are} ~ H_1{\rm {}-invariant} \\
X_1 &=& \{3,6,11,12,13,19,22,23,24,25,26,27,29,30,33,37,38,39,
41,43,44,45, \\
&& ~ 46,48,50,52,53,54,57,58\} \\
X_2 &=& \{3,11,12,13,14,15,19,22,25,26,31,35,37,38,41,43,44,46,
48,50,51, \\
&& ~ 52,55,56,58,60,61\} \\
X'_4 &=& 1,4,7,9,14,16,18,21,22,25,26,27,28\} \\
\\
&& (63;29,31,31,24;52) ~ (s{**}) \\
X'_1 &=& \{2,3,5,8,9,10,12,15,16,22,24,26,28,31\} \\
X_2 &=& \{0,1,2,3,4,6,7,8,10,11,12,15,20,23,24,28,29,30,31,34,
40,42,43,45, \\
&& ~ 49,50,52,58,59,60,61\} \\
X_4 &=& \{0,2,8,9,10,12,15,16,20,25,26,30,33,34,37,39,45,46,48,
50,57,60, \\
&& ~ 61,62\} \\
\\
&& (63;29,31,31,24;52) ~ ({**}s) \\
X_1 &=& \{1,2,3,11,15,18,19,21,22,26,27,28,29,30,33,34,35,36,38,
39,45,48, \\
&& ~ 49,51,52,55,59,60,61\} \\
X_2 &=& \{0,2,3,5,10,11,17,18,19,21,23,24,28,30,32,36,37,38,39,
40,41,42, \\
&& ~ 43,45,46,48,51,52,53,56,62\} \\
X'_4 &=& \{2,6,10,11,12,14,17,20,22,25,27,29\} \\
\\
&& (63;27,31,31,25;51) ~ (s{**}), ~
{\rm all} ~ X_i ~ {\rm are} ~ H_1{\rm {}-invariant} \\
X'_1 &=& \{2,3,7,8,10,12,14,15,21,23,28,29,31\} \\
X_2 &=& \{3,6,7,11,12,13,14,19,22,23,24,25,26,28,29,33,35,37,38,
41,42,43, \\
&& ~ 44,46,48,49,50,52,53,56,58\} \\
X_4 &=& \{3,5,12,13,15,17,19,20,22,23,25,26,29,30,37,38,39,41,
42,48,51,52, \\
&& ~ 53,57,60\} \\
\\
&& (63;27,31,31,25;51) ~ ({**}s) \\
X_1 &=& \{0,1,2,3,4,6,7,10,11,12,15,19,23,25,28,31,32,39,40,47,
49,51,52,53, \\
&& ~ 58,59,62\} \\
X_2 &=& \{0,1,6,7,9,10,13,14,16,18,19,23,24,28,30,35,37,38,40,41,
43,44,46, \\
&& ~ 47,48,49,50,54,59,61,62\} \\
X'_4 &=& \{1,2,4,6,10,14,16,17,19,21,27,28\} \\
\end{eqnarray*}

\begin{eqnarray*}
&& (63;27,29,29,26;48) ~ ({s{**}}), ~ 
{\rm all} ~ X_i ~ {\rm are} ~ H_2{\rm {}-invariant} \\
X'_1 &=& \{1,2,3,5,10,12,13,15,19,21,25,29,31\} \\
X_2 &=& \{2,7,9,10,13,14,18,20,21,26,27,28,29,32,35,36,40,42,
44,45,49,50, \\
&& ~ 52,53,54,55,56,59,61\} \\
X_4 &=& \{4,7,8,9,11,16,17,18,20,21,22,23,26,28,29,32,36,37,41,
42,43,44, \\
&& ~ 46,47,49,59\} \\
\\
&& (63;27,29,29,26;48) ~ ({**}s) \\
X_1 &=& \{0,1,3,5,6,8,9,12,15,17,21,25,26,27,28,29,31,35,36,41,
43,46,48,53, \\
&& ~ 54,57,61\} \\
X_2 &=& \{0,2,3,5,9,10,21,25,26,27,28,30,31,33,34,35,41,42,
43,44,45,46,49, \\
&& ~ 53,54,55,57,59,60\} \\
X'_4 &=& \{1,3,4,10,11,14,16,18,20,23,26,27,31\} \\
\\
\end{eqnarray*}

In the case $v=57$ our list contains two PDFs having different parameter sets and sharing the same symmetric block. The same 
is true for $v=61$.

\section{Acknowledgements}
This research was enabled in part by support provided by SHARCNET (http:// \\ 
www.sharcnet.ca) and Compute Canada (http://www.computecanada.ca).

\end{document}